\documentclass{article}

\RequirePackage[hyphens]{url}
\def\doi#1{\url{http://dx.doi.org/#1}}

\RequirePackage{amsmath}
\RequirePackage{amsthm}
\RequirePackage{amsfonts}
\RequirePackage{amssymb}

\newif\iftrimmarks  \trimmarksfalse \trimmarkstrue

\theoremstyle{theorem}

\theoremstyle{definition}

\begin{document}

\title{A SHORT PROOF OF A CONCRETE SUM}

\author{Samuel G. Moreno, Esther M. Garc\'{\i}a--Caballero\\{\small{\it Departamento de Matem\'aticas,
Universidad de Ja\'en}}\\{\small{\it 23071 Ja\'en, Spain}}}

\date{}

\maketitle

\begin{abstract}

We give an alternative proof of a formula that generalizes Hermite's identity. Instead involving modular arithmetic, our short proof relies on the Fourier-type expansion for the floor function and on a trigonometric formula.

\vspace{.5cm}



\end{abstract}

\vspace{.5cm}

 A closed form for $\sum_{k=0}^{m-1}\left\lfloor \frac{x+ n k}{m} \right\rfloor$, where  $x$ is a real number and $m$, $n$  are integers with $m>0$ can be found in \cite{GKP}, where the authors use modular arithmetic to establish the formula in an elementary (although somewhat long) manner. Our aim is to give a short proof of the above-mentioned result. To this end, define the {\it fractional part} of any real $x$ by $\{x\}=x-\lfloor x \rfloor$ and notice that it is a periodic piecewise linear function, discontinuous at each integer point, whose Fourier expansion gives us
\begin{eqnarray}\label{e005}
\lfloor x \rfloor=x-\frac{1}{2}+\frac{1}{\pi}\sum_{j=1}^{\infty}\frac{\sin (2 \pi j x)}{j},\qquad x\in \mathbb{R}\setminus \mathbb{Z}.
\end{eqnarray}
\noindent If $f(x)$ stands for the right-hand side of (\ref{e005}), then $f(n)=n-1/2$ for each integer $n$,  and thus $\lfloor n \rfloor=n=f(n)+1/2$. From (\ref{e005}) and using that $\sum_{k=0}^{m-1}k=(m-1)m/2$, one gets
\begin{eqnarray}\label{e010}
\sum_{k=0}^{m-1}\left\lfloor \frac{x+n k}{m} \right\rfloor=x+\frac{(m-1)n}{2}-\frac{m}{2}+\frac{1}{\pi}\sum_{j=1}^{\infty}\frac{\sum_{k=0}^{m-1}\sin \left(2 \pi j \left(\frac{x+n k}{m}\right)\right)}{j},
\end{eqnarray}

\noindent provided none of the $(x+ n k)/m$ is an integer. By using $\sum_{k=0}^p \sin (z+ ak)=\csc(a/2) \sin (a(p+1)/2)\sin (z+a p/2)$ (see \cite{FS}) we establish that
\begin{eqnarray}\label{e020}
\sum_{k=0}^{m-1}\sin \left(2 \pi j \left(\frac{x+n k}{m}\right)\right)=\frac{\sin (\pi j n)}{\sin (\pi j \frac{n}{m})}\sin (\pi j n-\pi j \frac{n}{m}+2 \pi j \frac{x}{m}).
\end{eqnarray}
\noindent Note that the above sum vanishes except, eventually, when the denominator at the right-hand side also vanishes. Therefore, denoting $d={\rm gcd}(m,n)$, $m'=m/d$ and $n'=n/d$, (\ref{e020}) may differ from zero only when $jn/m=jn'/m'\in \mathbb{Z}$, namely, when $j=lm'=lm/d$ ($l=1,2,3,\ldots$). With this in mind, (\ref{e010}) finally transforms to
\begin{eqnarray}\label{e025}
\sum_{k=0}^{m-1}\left\lfloor \frac{x+n k}{m} \right\rfloor&=&x+\frac{(m-1)n}{2}-\frac{m}{2}+\frac{d}{\pi}\sum_{l=1}^{\infty}\frac{\sum_{k=0}^{m-1}\sin \left(2 \pi l n' k+2 \pi l\displaystyle{\frac{x}{d}}\right)}{l m} \nonumber\\
&=&\frac{(m-1)(n-1)}{2}-\frac{1}{2}+\frac{d}{2}+d\left(\frac{x}{d}-\frac{1}{2}+\frac{1}{\pi}\sum_{l=1}^{\infty}\frac{\sin \left(2 \pi l\displaystyle{\frac{x}{d}}\right)}{l}\right)\nonumber\\
&=&\frac{(m-1)(n-1)}{2}+\frac{d-1}{2}+d\left\lfloor\frac{x}{d} \right\rfloor.
\end{eqnarray}

A final comment is in order. If $(x+ n k_0)/m\in \mathbb{Z}$ for some $0\leq k_0\leq m-1$, then it is readily verified that:
\begin{enumerate}
  \item $x/d\in \mathbb{Z}$. Effectively, if some $l\in \mathbb{Z}$ exists such that $(x+ n k_0)/m=l$, then  $(x+ n'd k_0)/(m'd)=l$, which implies that
  \begin{eqnarray*}
  \frac{x}{d}=m'l-n'k_0\in \mathbb{Z}.
  \end{eqnarray*}
  \item if $k_1\neq k_0$ verifies that $0\leq k_1\leq m-1$ and also $(x+ n k_1)/m=l_1\in \mathbb{Z}$, then $|k_2-k_1|$ is a multiple of $m'$. To check it, just observe that from $(x+ n k_0)/m=l_0$ and $(x+ n k_1)/m=l_1$ one gets $n (k_1-k_0)=m(l_1-l_0)$ or $n' (k_1-k_0)=m'(l_1-l_0)$; since ${\rm gcd}(m',n')=1$, then $n'$ divides $(l_1-l_0)$, so finally $(k_1-k_0)=s m'$ for some integer $s$.
  \item for each integer $r$ such that $k_r=k_0+r m'\in \{0,1,\ldots m-1\}$, it also holds $(x+ n k_r)/m\in \mathbb{Z}$. To show it, use
  \begin{eqnarray*}
  \frac{x+ n k_0}{m}=l_0, \qquad\mbox{ and } \qquad\frac{x+ n' d k_0}{m' d}=l_0,
  \end{eqnarray*}
  \noindent  to obtain
  \begin{eqnarray*}
 \frac{x}{d}&=&m'l_0-n'k_0=m'l_0+m'rn'-n'k_0-m' r n'\\&=&m'(l_0+r n')-n'(k_0+r m').
  \end{eqnarray*}
\end{enumerate}

\noindent These three items above show that there are exactly $d$ distinct $k_j$s in $\{0,1,\ldots,m-1\}$ for which $(x+ n k_j)/m\in \mathbb{Z}$. Thus, (\ref{e025}) holds true in this case too, because
\begin{eqnarray*}
\sum_{k=0}^{m-1}\left\lfloor \frac{x+n k}{m} \right\rfloor&=&\left(x+\frac{(m-1)n}{2}-\frac{m}{2}+\frac{d}{\pi}\sum_{l=1}^{\infty}\frac{\sum_{k=0}^{m-1}\sin \left(2 \pi l n' k+2 \pi l\displaystyle{\frac{x}{d}}\right)}{l m}\right)\\&+&\frac{d}{2}
\end{eqnarray*}
\noindent (we have added $d/2$ to correct the value of $f(\cdot)$ in the $d$ cases in which its argument is an  integer); therefore
\begin{eqnarray*}
\sum_{k=0}^{m-1}\left\lfloor \frac{x+n k}{m} \right\rfloor&=&\frac{(m-1)(n-1)}{2}-\frac{1}{2}+\frac{d}{2}\\&+&d\left(\frac{x}{d}-\frac{1}{2}+\frac{1}{\pi}\sum_{l=1}^{\infty}\frac{\sin \left(2 \pi l\displaystyle{\frac{x}{d}}\right)}{l}+\frac{1}{2}\right)\\&=&\frac{(m-1)(n-1)}{2}+\frac{d-1}{2}+d\left(f\left(\frac{x}{d}\right)+\frac{1}{2}\right)\\
&=&\frac{(m-1)(n-1)}{2}+\frac{d-1}{2}+d\left\lfloor\frac{x}{d} \right\rfloor,
\end{eqnarray*}

\noindent and we are done.

\vfill\eject

\end{document}